\documentclass[conference]{IEEEtran}
\usepackage{subfigure}
\usepackage{lineno,hyperref}
\usepackage[table,xcdraw]{xcolor}
\usepackage{amsmath}
\usepackage{cases}
\usepackage{graphicx}
\usepackage[margin=1in]{geometry}
\usepackage{color}
\usepackage{array}
\usepackage{multirow}
\newcolumntype{P}[1]{>{\centering\arraybackslash}p{#1}}
\usepackage{caption}% http://ctan.org/pkg/caption
\begin{document}
\title{A note on the Bures-Wasserstein metric}
\author{Shravan Mohan\\ 17-004, Mantri Residency, Bannerghatta Main Road, Bangalore.  
}
\newgeometry{top=1in,bottom=0.75in,right=0.75in,left=0.75in}
\maketitle

\begin{abstract}
In this brief note, it is shown that the Bures-Wasserstein (BW) metric on the space positive definite matrices lends itself to convex optimization. In other words, the computation of the BW metric can be posed as a convex optimization problem. In turn, this leads to efficient computations of (i) the BW distance between convex subsets of positive definite matrices, (ii) the BW barycenter, and (iii) incorporating BW distance from a given matrix as a convex constraint. Computations are provided for corroboration. 
\end{abstract}
\section{Introduction}
\noindent Consider the set of positive definite matrices of dimension $n$ given by $\mathcal{P}(n)$. The Bures-Wasserstein (BW) metric between $A$ and $B$ in $\mathcal{P}(n)$ is given by the closed form \cite{bhatia2019bures}:
\begin{align}
    \rho^2(A, B) = \mbox{Tr}(A) + \mbox{Tr}(B) - 2\mbox{Tr}\left(\sqrt{\sqrt{A}B\sqrt{A}}\right).
\end{align}
Here, $\sqrt{X}$ denotes the unique symmetric square root of a positive definite matrix $X$. That is:
\begin{align}
    \sqrt{X} = U\Sigma^{\frac{1}{2}}U^\top,
\end{align}
where $X = U\Sigma U^\top$ is the singular value decomposition of $X$ and $\Sigma^{\frac{1}{2}}$ is the 
element-wise square root of $\Sigma$.\\\\
The BW metric is defined in the following way. The set of square roots of a positive definite matrix $X$ is given by:
\begin{align}
    \left\{ \sqrt{X}U: U \in O(n)\right\},
\end{align}
where $O(n)$ is the set of real unitary matrices of dimension $n$. Then, for two positive definite matrices $A$ and $B$, the Frobenius distance ($||.||_F$) between the sets of their respective square roots is defined as the BW metric. Mathematically, this gives:
\begin{align}
    \rho(A, B) = \min_{U, V\in O(n)}\left|\left|\sqrt{A}U - \sqrt{B}V\right|\right|_F
\end{align}
\section{BW metric \& Semidefinite Programming}
Since the Frobenius norm is unitary invariant, the BW metric can also be written as:
\begin{align}
    \rho(A, B) = \min_{V\in O(n)}\left|\left|\sqrt{A} - \sqrt{B}V\right|\right|_F,
\end{align}
Thus,
\begin{align}
    \rho^2(A, B) = \min_{U\in O(n)}\mbox{Tr}(A) + \mbox{Tr}(B) - 2\mbox{Tr}\left(\sqrt{A}\sqrt{B}U\right).
\end{align}
Now, the following well-known lemma comes to the aid for solving the above optimization problem as a convex optimization problem \cite{boyd2004convex}.\\
\textbf{Lemma}: The linear SDP given by
\begin{align}
    &\max_{U}~~~~\mbox{Tr}(KU)\\
    &\mbox{subject to~~} \begin{bmatrix}
        I & U^\top\\
        U & I
    \end{bmatrix} \succeq 0. 
\end{align}
has an optimal solution $\tilde{U}$ such that $\tilde{U}^\top \tilde{U} = I$. \\
\textbf{\textit{Proof}}: Firstly, note that if a feasible point $U$ is such that some diagonal elements of $KU$ are negative, then the matrix $UD$, where  $D$ is diagonal such that
$$
D_{i,i}= 
\begin{cases}
    1,& \text{if } (KU)_{i,i}\geq 0\\
    -1,              & \text{otherwise},
\end{cases}
$$
would yield a higher cost function value. Thus, at optimality, the diagonal elements of $KU$ are non-negative. Also note that $DU$ satisfies the semidefinite constraint if $U$ does. Secondly, suppose the optima $\tilde{U}$ was such that $I \succ \tilde{U}^\top \tilde{U}$. Let $\tilde{U} = PSQ^\top$ by SVD. By our assumption, some of the elements of $S$ have to be zero or less than 1. Now consider the unitary matrix $U = PD Q^\top $, where $D$ is a diagonal matrix (from the first observation) which makes all the diagonal elements of $KPD Q^\top$ non-negative. 
Also note that
\begin{align}\small
\mbox{Tr}(KPDQ^\top )  = \mbox{Tr}(Q^\top KPD)\geq \mbox{Tr} (Q^\top KPS) \nonumber = \mbox{Tr}\left(KU\right)\nonumber,
\end{align}
which finally implies that the convex relaxation is tight. The result also applies to the case where the constraint is:  
\begin{align}
    \begin{bmatrix}
        G & U^\top \\
        U & I
    \end{bmatrix} \succeq 0,
\end{align}
since this is equivalent to the constraint 
\begin{align}
    \begin{bmatrix}
        I & \sqrt{G^{-1}}U^\top \\
        U\sqrt{G^{-1}} & I
    \end{bmatrix} \succeq 0,
\end{align}
With the above lemma, the computation of BW distance can also be written as:
\begin{align}
    \rho^2(A, B) =~~& \min~~ \mbox{Tr}(A) + \mbox{Tr}(B) - 2\mbox{Tr}\left(\sqrt{A}U\right)\\
    &\mbox{subject to~~} \begin{bmatrix}
    B & U^\top\\
    U & I
    \end{bmatrix}\succeq 0.
\end{align}
\begin{figure*}[t!]
\centering
\fbox{\begin{minipage}{7in}
\begin{enumerate}
    \item{
    Computation of BW distance between convex subsets $\mathcal{A}$ and $\mathcal{B}$ of PD matrices:
    \begin{itemize}
        \item Choose any matrix $A\in \mathcal{A}$.
        \item{Compute the following:
            \begin{align}
                B =& \mbox{~arg}\min_{X} \mbox{Tr}(X) - 2\mbox{Tr}\left(\sqrt{A}K\right)\\
                &\mbox{subject to~} \begin{bmatrix}
                    X & K^\top\\
                    K & I
                \end{bmatrix} \succeq 0,~ X \in \mathcal{B}.
            \end{align}    
        }
        \item{Compute the following:
            \begin{align}
                A =& \mbox{~arg}\min_{X} \mbox{Tr}(X) - 2\mbox{Tr}\left(\sqrt{B}K\right)\\
                &\mbox{subject to~} \begin{bmatrix}
                    X & K^\top\\
                    K & I
                \end{bmatrix} \succeq 0,~ X \in \mathcal{A}.
            \end{align}    
        }
        \item Repeat the above two steps till convergence. 
    \end{itemize}
    }
    \item {
    Computation of the BW barycenter:
    \begin{align}
    &\min_{X}~~ \sum_{i=1}^N w_i \left( \mbox{Tr}(A_i) + \mbox{Tr}(X) - 2\mbox{Tr}( \sqrt{A_i}K_i)\right)\\
    &\mbox{subject to~~} \begin{bmatrix}
        X & K_i^\top\\
        K_i & I
    \end{bmatrix} \succeq 0, ~ \forall i.
    \end{align}
    }
    \item{
    Incorporating BW distance from a matrix as a constraint in a convex program:
    \begin{align}
        &\min_{X} ~~~~f\left(X\right); ~~f \mbox{~convex in $X$}\\
        &\mbox{subject to~~}
        X\in \mathcal{C}, ~~\begin{bmatrix}
        X & K^\top\\
        K & I
        \end{bmatrix} \succeq 0, \\
        &~~~~~~\mbox{Tr}(A) + \mbox{Tr}(X) - 2\mbox{Tr}\left(\sqrt{A_i}K\right) \leq d_{A_i}, \forall i. 
    \end{align}
    }
\end{enumerate}
\end{minipage}}
\caption{The Main Algorithms}
\label{fig:minmaxalgo}
\end{figure*}

\begin{table*}[h!]
\centering
\caption{Parameters and computational results.}
\resizebox{6.9in}{!}{\begin{tabular}{|c|c|c|}
\hline
Purpose & Parameters & Results   \\ \hline
BW distance between convex subsets of PD matrices & $\mathcal{A}=\left\{X\in S_5^+|\mbox{Tr}(X)=1\right\}$, $\mathcal{B}=\left\{X\in S_5^+|\mbox{Tr}(X)=2\right\}$ & \textcolor{blue}{$A = \begin{bmatrix}
         0.3209& -0.1364& -0.1069& -0.1686&  0.0726\\
       -0.1364&  0.5256&  0.1634& -0.0637& -0.1171\\
       -0.1069&  0.1634&  0.5295&  0.0262& -0.095 \\
       -0.1686& -0.0637&  0.0262&  0.2931&  0.0048\\
        0.0726& -0.1171& -0.095 &  0.0048&  0.3308
\end{bmatrix}$,  $B = \begin{bmatrix}
     0.1605& -0.0682& -0.0535& -0.0843&  0.0363\\
       -0.0682&  0.2628&  0.0817& -0.0319& -0.0585\\
       -0.0535&  0.0817&  0.2647&  0.0131& -0.0475\\
       -0.0843& -0.0319&  0.0131&  0.1466&  0.0024\\
        0.0363& -0.0585& -0.0475&  0.0024&  0.1654
\end{bmatrix}$} \\ \hline
 & $w = \left[0.8766, 0.6682, 1.0852, 1.1009, 0.524\right]$ & \\
 & $A1 = \begin{bmatrix}
     2.7273& -1.3426& -1.4873&  1.1069& -0.5844\\
    -1.3426&  5.6047& -0.7192&  0.3519&  1.0648\\
       -1.4873& -0.7192&  4.6821& -0.9547& -1.6117\\
        1.1069&  0.3519& -0.9547&  2.4089& -0.9744\\
       -0.5844&  1.0648& -1.6117& -0.9744&  3.5771 
\end{bmatrix}$, $A2 = \begin{bmatrix}
     5.6143& -0.1039&  1.4161&  0.0105&  0.9256\\
       -0.1039&  4.6277&  0.5304&  0.0571&  0.6138\\
        1.4161&  0.5304&  7.017 & -0.4625& -0.3483\\
        0.0105&  0.0571& -0.4625&  5.4935&  1.2015\\
        0.9256&  0.6138& -0.3483&  1.2015&  8.2474 
\end{bmatrix}$ & 
\textcolor{blue}{$\textcolor{blue}{X_{\mbox{opt}} = \begin{bmatrix}
    3.8514& -0.5993& -0.0722&  0.5644& -0.4899\\
       -0.5993&  4.8924&  0.1755&  0.0716&  0.2198\\
       -0.0722&  0.1755&  4.4109& -0.1818& -0.6419\\
        0.5644&  0.0716& -0.1818&  3.9922& -0.3331\\
       -0.4899&  0.2198& -0.6419& -0.3331&  4.5659 
\end{bmatrix}}$}\\
BW barycenter  & $A3 = \begin{bmatrix}
    5.4601& -0.1268& -0.7682& -0.729 & -0.909 \\
       -0.1268&  7.7425&  0.1735&  0.4499& -0.511 \\
       -0.7682&  0.1735&  6.8627& -0.3396& -1.259 \\
       -0.729 &  0.4499& -0.3396&  6.7328&  1.1921\\
       -0.909 & -0.511 & -1.259 &  1.1921&  4.2019 
\end{bmatrix}$, $A4 = \begin{bmatrix}
    2.937 & -1.1282&  0.3996&  0.9282& -0.3372\\
       -1.1282&  3.3586& -0.4808& -1.112 &  0.3812\\
        0.3996& -0.4808&  2.1708&  0.4026& -0.1732\\
        0.9282& -1.112 &  0.4026&  3.043 & -0.8748\\
       -0.3372&  0.3812& -0.1732& -0.8748&  4.4907 
\end{bmatrix}$ & $X_{\mbox{opt}}$ is computed using the optimization algorithm, while $X_{\mbox{fp}}$ is calculated using the fixed point equation \\
 & $A5 = \begin{bmatrix}
    4.5401&  1.2074&  1.3077&  1.6847& -1.2072\\
        1.2074&  3.9336&  2.5037&  1.5876& -0.3888\\
        1.3077&  2.5037&  3.8015&  0.5648&  0.9108\\
        1.6847&  1.5876&  0.5648&  4.1194& -1.8946\\
       -1.2072& -0.3888&  0.9108& -1.8946&  4.6055 
\end{bmatrix}$ & 
$\textcolor{blue}{X_{\mbox{fp}} = \begin{bmatrix}
     3.8514& -0.5994& -0.0722&  0.5644& -0.4899\\
       -0.5994&  4.8924&  0.1756&  0.0716&  0.2198\\
       -0.0722&  0.1756&  4.4108& -0.1818& -0.6419\\
        0.5644&  0.0716& -0.1818&  3.9921& -0.3331\\
       -0.4899&  0.2198& -0.6419& -0.3331&  4.5658 
\end{bmatrix}}$
\\ \hline
BW distance as convex constraint &
$f(X) = \left|\left|X\right|\right|_{F}$,~ $\rho^2(A,X)\leq 10$, ~$A = \begin{bmatrix}
      6.5722& -0.4557&  0.018 &  0.0854&  0.1883\\
       -0.4557&  6.3399& -0.0739& -0.1726& -0.2416\\
        0.018 & -0.0739&  5.8477& -0.2659& -0.2295\\
        0.0854& -0.1726& -0.2659&  5.5408& -0.3855\\
        0.1883& -0.2416& -0.2295& -0.3855&  5.6995 
\end{bmatrix}$
& \textcolor{blue}{$X = \begin{bmatrix}
    1.1203& -0.036 &  0.0016&  0.0071&  0.0152\\
    -0.036 &  1.1016& -0.0065& -0.0149& -0.0201\\
    0.0016& -0.0065&  1.0617& -0.0234& -0.0202\\
    0.0071& -0.0149& -0.0234&  1.0346& -0.0341\\
    0.0152& -0.0201& -0.0202& -0.0341&  1.0482 
\end{bmatrix}$}\\ \hline 
\end{tabular}}
\label{tab:min_max}
\end{table*}
\noindent \textbf{\underline{BW distance between convex sets}}: The convex routine to compute the BW metric can be used to find the BW distance between two convex subsets of positive definite matrices using alternating projections \cite{bauschke1996projection}. The alternating projections method proceeds by finding the distance of a point from one set to the other  alternatively, while latching onto the latest iterate. Recall that this method converges to the distance between convex sets (and the corresponding matrices in the two subsets), given a metric on the point set. More precisely, consider two convex sets of positive definite matrices $\mathcal{A}$ and $\mathcal{B}$. Then, the alternating steps would be:
\begin{align}
    &\min_{K, A}~~ \mbox{Tr}(A) + \mbox{Tr}(B) - 2\mbox{Tr}(\sqrt{B}K)\\
    &\mbox{subject to~~}     \begin{bmatrix}
        A & K^\top\\
        K & I
    \end{bmatrix} \succeq 0, ~ A\in \mathcal{A}.
\end{align}
and 
\begin{align}
    &\min_{K, B}~~ \mbox{Tr}(A) + \mbox{Tr}(B) - 2\mbox{Tr}( \sqrt{A}K)\\
    &\mbox{subject to~~}     \begin{bmatrix}
        B & K^\top\\
        K & I
    \end{bmatrix} \succeq 0, ~ A\in \mathcal{A}.
\end{align}
\textbf{\underline{BW barycenter}}: The same idea can be applied towards finding a weighted BW barycenter. Conventionally, the BW barycenter is computed as a solution to a fixed point equation given by:
\begin{align}
    X = \sqrt{X^{-1}}\left(\sum_{i=1}^N w_i \sqrt{\sqrt{X} A_i \sqrt{X}}\right)^2\sqrt{X^{-1}}.
\end{align}
Although, the fixed point iteration converges to the BW barycenter, a convex approach is easy to understand and allows inclusion of convex constraints on the barycenter. Consider the convex optimization problem:
\begin{align}
    &\min_{X}~~ \sum_{i=1}^N w_i \left( \mbox{Tr}(A_i) + \mbox{Tr}(X) - 2\mbox{Tr}( \sqrt{A_i}K_i)\right)\\
    &\mbox{subject to~~}  \begin{bmatrix}
        X & K_i^\top\\
        K_i & I
    \end{bmatrix} \succeq 0, ~ \forall i.
\end{align}
Note that the optima in this case too has to lie on the boundary of each constraint (by the aforementioned lemma), and hence the BW barycenter can be calculated this way. \\\\
\textbf{\underline{BW distance as convex constraint}}: Consider the constraint set given by:
\begin{align}
    \left\{X\in \mbox{PD}(n) ~~|~~ \rho(A,X) \leq d \right\}.
\end{align}
This set can be represented as a convex constraint:
\begin{align}
    \left\{X\in PD(n), ~K\in R^{n,n} ~~|~~ 
    \begin{bmatrix}
        X & K^\top \\
        K & I
    \end{bmatrix}\succeq 0 \right. ~\&~\\
    \left. \mbox{Tr}(A) + \mbox{Tr}(X) -2\mbox{Tr}(\sqrt{A} K)\leq d^2 
    \right\}.
\end{align}
Note that if $(X, K)$ belongs to the constraint set, then obviously $$
\min_{X, K\in O(n)} ~\left(\mbox{Tr}(A) + \mbox{Tr}(X) -2\mbox{Tr}(\sqrt{A} K)\right) = \rho^2(A,X) \leq d^2.
$$
An example of convex optimization problem incorporating such a constraint would be:
\begin{align}
    &\min_X ~~~~||X||_{{F}}\\
    &\mbox{sub to~~} \rho(A,X) \leq d.
\end{align}
\section{Computations}
Computational results for the three use cases are shown in Table 1. For the first use case, the convex subsets of positive definite matrices are BW distance between two convex subsets of positive definite matrices are $\mathcal{A}=\left\{X\in P(n)|\mbox{Tr}(X)=1\right\}$ and  $\mathcal{B}=\left\{X\in P(n)|\mbox{Tr}(X)=2\right\}$. For the second use case, the weight vector and the matrices of which the BW barycenter needs to be calculated are presented. The computation of the barycenter using the fixed point equation also yields essentially the same result, thereby corroborating this paper's claim. For the third case, the example optimization problem is chosen with the matrix $A$ given in the table and $d$ set as $\sqrt{10}$. All the computations were done using CVXPY \cite{diamond2016cvxpy}.
\section{Conclusion}
In this paper it was shown that the computation of the BW metric can be done using convex optimization. This resulted in numerically efficient routines for calculating the BW distance between convex subsets of matrices, the BW barycenter of a finite set of positive definite matrices and incorporating BW distance from a matrix as a convex constraint in a convex optimization problem. Computational examples were provided for corroboration.

% \begin{figure*}[t!]
% \centering
% \fbox{\begin{minipage}{7in}
% \begin{enumerate}
%     \item {For singleton sets, 2-tuple \& 3-tuple combinations of circles, find the best possible packing. For each combination:
%     \begin{enumerate}
%         \item Generate a fine grid and a circle intersection graph. 
%         \item Score (heuristic) nodes based on proximity to the corners / edges / 
%         to the already placed circles, and the increase in packing density. 
%         \item Choose a node with the highest score and remove its neighbours. 
%         \item Repeat last two steps till no nodes are left, calculate the number of cylinders than can be stacked in every slot of the rack. \item Store those numbers in matrices $G_1$, $G_2$ and $G_3$; the number implying the tuple size. The row id corresponds to the exact tuple, while the column tuple denotes the slot number. 
%     \end{enumerate}
%     }
%     \item {Define $E_1$, $E_2$ and $E_3$ as before. Perform the following optimization.
%     \begin{align}
%        & \min ~~~\max_{Z_1, Z_2, Z_3} \left\{\boldsymbol{1}^\top Z_1 + \boldsymbol{1}^\top Z_2 + \boldsymbol{1}^\top Z_3\right\}\\
%        & \mbox{subject to}\\
%        & ~~~~~~~~\left(E_1G_1 + E_2G_2 + E_3G_3\right)\boldsymbol{1} \geq V\\
%        & ~~~~~~~~~~~~~~ 0\leq Z_1, Z_2, Z_3 \leq 1.
%     \end{align}
%     } 
%     \item Round off $Z_1, Z_2, Z_3$ to their celing. 
% \end{enumerate}
% \end{minipage}}
% \caption{The Main Algorithms}
% \label{fig:minmaxalgo}
% \end{figure*}

\bibliographystyle{IEEEtran}
\bibliography{biblio}
\end{document}